\documentclass[10pt,twoside]{amsart}
\usepackage{amsmath, amsthm, amscd, amsfonts, amssymb, graphicx}%, color}
\numberwithin{equation}{section}

\begin{document}
\begin{center}\large{{\bf{Lacunary ideal convergence in probabilistic normed spaces}}}

\vskip 0.5 cm
%\author
Bipan Hazarika and Ayhan Esi$^{\ast}$\\
%}}%{\footnote {The corresponding author}}}
Department of Mathematics, Rajiv Gandhi University, Rono Hills, Doimukh-791112. Arunachal Pradesh, India\\

Adiyaman University, Science and Art Faculty, Department of Mathematics, 02040, Adiyaman, Turkey\\

E-mail: bh\_rgu$@$yahoo.co.in; aesi23@hotmail.com.
\end{center}
\title{}
\author{}
\thanks{{May 13, 2014.
$^{\ast}$ Corresponding author.}}

%\title[Lacunary ideal convergence ...]{Lacunary ideal convergence in probabilistic normed spaces}
%\author{Bipan Hazarika}
%\address{%Adiyaman University, Science and Art Faculty, Department of Mathematics, 02040, Adiyaman, Turkey\\ 
%Department of Mathematics, Rajiv Gandhi University, Rono Hills, Doimukh-791112. Arunachal Pradesh, India\\}
%\email{ bh\_rgu@yahoo.co.in }
%\subjclass[2000]{40G15; 46S70; 54E70}
%\keywords{Ideal convergence; probabilistic normed space; lacunary sequence; $\theta$-convergence.}
%\date{{May 13, 2012}}

\begin{abstract}

An ideal $I$ is a family of subsets of positive integers $\mathbb{N}$ which is closed under taking finite unions and subsets of its elements. A sequence $(x_k)$ of real numbers is said to be lacunary $I$-convergent to a real number $\ell$, if for each $ \varepsilon> 0$ the set $$\left\{r\in \mathbb{N}:\frac{1}{h_r}\sum_{k\in J_r} |x_{k}-\ell|\geq \varepsilon\right\}$$ belongs to $I.$ The aim of this paper is to study the notion of lacunary $I$-convergence in probabilistic normed spaces as a variant of the notion
of ideal convergence. Also lacunary $I$-limit points and lacunary $I$-cluster points have been defined and the relation between them has been
established. Furthermore, lacunary-Cauchy and lacunary $I$-Cauchy sequences are
introduced and studied. Finally, we provided example which shows that our method of convergence in  probabilistic normed spaces is more general.\\

AMS Subjclass[2000]: 40G15; 46S70; 54E70\\


\end{abstract}
\maketitle
\pagestyle{myheadings}
\markboth{\rightline {\scriptsize  Hazarika, Esi}}
         {\leftline{\scriptsize Lacunary ideal convergence in...}}

\maketitle

\section{Introduction}
\normalfont
%The concept of ideal convergence as a generalization of statistical convergence, and any concept involving statistical convergence play a vital role not only in pure mathematics but also in other branches of science involving mathematics, especially in information theory, computer science, biological science, dynamical systems, geographic information systems, population modelling, and motion planning in robotics.\newline

Steinhaus \cite{Steinhaus} and
Fast \cite{Fast} independently introduced the notion of statistical convergence for sequences of real numbers. %, the theory has been investigated and developed by several authors (\cite{Connor},  \cite{Fridy85}, \cite{Fridy93},  \cite{Salat80}, \cite{SalatTijdeman}).  
Over the years  and under different names statistical convergence has been discussed in the theory of Fourier analysis, ergodic theory  and number theory. Later on it was further investigated from various points of view. For example, statistical convergence has been investigated in summability theory by ( Connor \cite{Connor}, Fridy \cite{Fridy85},  \u{S}al\'at \cite{Salat80}), number theory and mathematical analysis by  (Buck \cite{Buck}, Mitrinovi\'c et al., \cite{MitrinovicSandorCrstici}), 
topological groups (\c{C}akalli (\cite{Cakalli96, Cakalli09})), topological spaces (Di Maio and Ko\u{c}inac \cite{MaioKocinac08}), function spaces (Caserta and Ko\u{c}inac \cite{CasertaMaioKocinac}), locally convex spaces (Maddox\cite{Maddox88}), measure theory (Cheng et al., \cite{ChengLinLanLiu}, Connor and Swardson \cite{Con}, Miller\cite{Miller95}). Fridy and Orhan \cite{FridyOrhan93} introduced the concept of lacunary statistical convergence. Some work on lacunary statistical convergence can be found in (\cite{Cakalli96, FridyOrhanLacunary statistical summability, HazarikaSavas123,  Li00}). \\
%In the recent years, generalization of statistical  convergence have appeared  in the study of strong integral summability and the structure of ideals of bounded continuous  functions \cite{Con}. 
%Mursaleen \cite{Mursaleen}, introduced the $\lambda$-statistical convergence for real sequences. In this article, we consider only sequences of real numbers, so  that "a sequence" means "a sequence of real numbers". 

Kostyrko, et. al \cite{KostyrkoSalatWilczynski} introduced the notion of $I$-convergence as a generalization of statistical convergence which is based on the structure of an admissible ideal $I$ of subset of natural numbers $\mathbb{N}$. Kostyrko, et. al \cite{KostyrkoMacajSalatSleziak} gave some of basic properties of $I$-convergence and dealt with extremal $I$-limit points. Further details on ideal convergence can be found in (\cite{CakalliHazarika, Dems, EsiHazarika, EsiHazarika00, HazarikaSavas, Hazarika622, Hazarika12, Hazarika769, HazarikaKimarGuillen, HazarikaOn ideal convergence in topological groups, LahiriDas, SalatTripathyZiman04, TripathyHazarika, TripathyHazarika09, TripathyHazarika11}), and many others. The notion of lacunary ideal convergence of real sequences was introduced in (\cite{Choudhary09, TripathyHazarikaChoudhary}) and Hazarika (\cite{Hazarika09, Hazarika}),  was introduced the lacunary ideal convergent sequences of fuzzy real numbers and studied some properties. Debnath \cite{Debnath} introduced the notion lacunary ideal convergence in intuitionistic fuzzy normed linear spaces. Recently, Yamanci and G\"{u}rdal \cite{YamanciGurdal} introduced the notion lacunary ideal convergence in random $n$-normed space. \\
%Das et. al  \cite{DasKostyrkoWilczynskiMalik} introduced the concept of $I$-convergence of double sequences in metric space and studied some properties of this convergence. A lot of developments have been made in this areas after the works of \cite{Kumar}, \cite{TripathyTripathy}. 

%Although an ideal is defined as a hereditary and additive family of subsets of a non-empty arbitrary set $X$, here in our study it suffices to take $I$ as
 A family of subsets of $\mathbb{N}$, positive integers, i.e. $I\subset2^{\mathbb{N}}$ is an ideal on $\mathbb{N}$ if and only if 
 \begin{enumerate}
 \item[(i)] $\phi\in I,$
 \item[(ii)] $A\cup B\in I$ for each $ A,B\in I,$ 
 \item[(iii)] each subset of an element of $I$ is an element of $I$. 
 \end{enumerate}
 A non-empty family of sets $F\subset2^{\mathbb{N}}$ is a filter on $\mathbb{N}$ if and only if 
 \begin{enumerate}
 \item[(a)] $\phi \notin F$
 \item[(b)] $A\cap B\in F$ for each $A,B\in F,$
 \item[(c)] any subset of an element of $F$ is in $F$.
 \end{enumerate}
  An ideal $I$ is called \textit{non-trivial} if $I\neq\phi$ and $\mathbb{N}\notin I.$ Clearly $I$ is a non-trivial ideal if and only if $F=F(I) =\{\mathbb{N}-A: A\in I\} $ is a filter in $\mathbb{N},$ called the filter associated with the
ideal $I.$\\

 A non-trivial ideal $I$ is called \textit{admissible} if and only if $\{ \{ n \} :n\in \mathbb{N} \}\subset I$. A non-trivial ideal $I$ is \textit{maximal} if there cannot exists any non-trivial
ideal $J \neq I$ containing $I$ as a subset.\\
 %Further details on ideals can be found in Kostyrko, et.al (see \cite{KostyrkoSalatWilczynski}).
%Throughout this paper we assume $I$ is a non-trivial admissible ideal in $\mathbf{N}.$ 
Recall that a sequence $x=( x_{k})$ of points in $\mathbb{R}$ is said to be $I$-convergent to a real number $\ell$ if $ \{k\in \mathbb{N}: \vert x_k -\ell\vert \geq\varepsilon \} \in \textit{I}$ for every $\varepsilon >0$ (\cite{KostyrkoSalatWilczynski}). In this case we write $ I-\lim x_k = \ell $.\\%The notion was further investigated by \u{S}al\'at, et.al \cite{SalatTripathyZiman04}} and others.\\

By a lacunary sequence $\theta = (k_r),$ where $k_0 = 0$ , we shall mean an
increasing sequence of non-negative integers with $k_r - k_{r-1}\rightarrow
\infty$ as $r\rightarrow \infty.$ The intervals determined by $\theta$ will
be denoted by $J_r= (k_{r-1}, k_r]$ and we let $h_r = k_r - k_{r-1}.$ The
space of lacunary strongly convergent sequences $\mathcal{N}_\theta$ was defined by
Freedman et al. \cite{FreedmanSemberRaphael78} as follows:\newline
\[
\mathcal{N}_\theta = \left\{ x = (x_k): \lim_{r} \frac{1}{h_r} \sum\limits_{k\in J_r}
\vert x_k - L \vert =0, ~\mbox{for ~some} ~L \right \}.
\]

Menger \cite{Menger} proposed the probabilistic concept of the distance by replacing the number $d(p, q)$ as the distance between points $p, q$ by a probability distribution
function $F_{p,q}(x)$. He interpreted $F_{p,q}(x)$ as the probability that the distance between $p$ and $q$
is less than $x.$ This led to the development of the area now called probabilistic metric spaces.
This is \u{S}herstnev \cite{Sherstnev} who first used this idea of Menger to introduce the concept of a PN
space. In 1993, Alsina et al. \cite{AlsinaSchweizerSklar93} presented a new definition of probabilistic normed space
which includes the definition of ¢Sherstnev as a special case. For an extensive view on this
subject, we refer (\cite{AlsinaSchweizerSklar97, ConstantinIstratescu, GuillenLallenaSempi, GuillenSempi, Karakus, KumarGuillen,  LafuerzaLallenaSempi, SchweizerSklar60, SchweizerSklar83}). %Quite recently, Karakus \cite{Karakus} defined statistical analogues of convergence and Cauchy double sequences on PN spaces and gave a useful characterization.
Subsequently, Mursaleen and Mohiuddine \cite{MursaleenMohiuddine} and Rahmat\cite{Rahmat} studied the ideal convergence in probabilistic normed spaces and V. Kumar and K. Kumar \cite{KumarKumar08} studied $I$-Cauchy and $I^{\ast}$-Cauchy sequences in probabilistic normed spaces.\\

%The notion of statistical convergence was introduced by Fast \cite{Fast51} and Schoenberg \cite{Schoenberg59} independently. A lot of developments have been made in this areas after the works of \u{S}al\'at \cite{Salat80} and Fridy \cite{Fridy85}. Over the years  and under different names statistical convergence has been discussed in the theory of Fourier analysis, ergodic theory  and number theory. 
%This concept was extended to the double sequences by Mursaleen and Edely \cite{MursaleenEdelyStatisticalconvergenceofdoublesequences} and Tripathy \cite{Tripathy} independently. 
%In the recent years, generalization of statistical  convergence have appeared  in the study of strong integral summability and the structure of ideals of bounded continuous  functions on Stone-\u{C}ech  compactification of the natural numbers. Moreover statistical convergence is closely related to the concept of convergence in probability, (see \cite{Con}). \newline

The notion of statistical convergence depends on the density (asymptotic or natural) of subsets of
$\mathbb{N}.$ A subset of $\mathbb{N}$ is said to have natural density $\delta\left(  E\right)  $ if
\[
\delta\left(  E\right)  =\lim_{n\rightarrow\infty}\frac{1}{n}|\{k\leq n: k\in E\}|   \mbox{~exists.}
\]
%\sum_{k=1}^{n}\chi_{E}\left(  k\right)
\textbf{Definition 1.1.} A sequence $x=\left(  x_{k}\right)  $ is said to be \textit{statistically
convergent} to $\ell$ if for every $\varepsilon>0$
\[
\delta\left(  \left\{  k\in\mathbb{N} :\left\vert x_{k}-\ell\right\vert \geq\varepsilon\right\}  \right) =0.
\]
In this case, we write $S-\lim x = \ell$ or $x_k \rightarrow \ell(S)$ and $S$ denotes the set of all statistically convergent sequences.\newline

\textbf{Definition 1.2.}(\cite{Choudhary09, TripathyHazarikaChoudhary}) Let $I\subset 2^{\mathbb{N}}$ be a non-trivial ideal. A real sequence
$x=\left(  x_{k}\right)  $  is said to be {\it lacunary I-convergent} or $I_{\theta}$-{\it convergent} to
$L\in \mathbb{R}$ if, for every $\varepsilon>0 $   the set%
\[
\left\{  r\in%
%TCIMACRO{\U{2115} }%
%BeginExpansion
\mathbb{N}
%EndExpansion
:\text{ }\frac{1}{h_{r}}\sum_{k\in J_{r}} |x_{k}-L| \geq \varepsilon  \right\}  \in I.
\]
$L$ is called the $I_{\theta}$-limit of the sequence $x=\left(  x_{k}\right),  $ 
and we write $I_{\theta}-\lim x=L.$\\

In this paper we study the concept of lacunary $I$-convergence in probabilistic normed spaces. We also define lacunary $I$-limit points and lacunary $I$-cluster points in probabilistic normed space and prove some interesting results.

%We recall some notations and basic definitions used in this paper.

%The existing literature on ideal convergence and its generalizations appears to have been restricted to real or complex sequences, but in recent years these ideas have been also extended to the sequences in fuzzy normed \cite{KumarKumar08} and intutionistic fuzzy normed spaces \cite{KumarKumar}, \cite{MursaleenMohiuddineEdely}, \cite{MohiuddineLohani}, \cite{KumarMursaleen}, \cite{MursaleenMohiuddineOnlacunarystatisticalconvergencewithrespecttotheintuitionisticfuzzynormedspace}, \cite{KarakusDemirciDuman}. Further details on ideal convergence can be found in \cite{Hazarika},  \cite{LahiriDas}, \cite{SahinerGurdalSaltanGunawan},  \cite{TripathyHazarika} and many others.\\

\section{Basic definitions and notations}

\qquad Now we recall some notations and basic definitions that we are going to
use in this paper.\\

%\textbf{Definition 1.2.}\cite{KostyrkoSalatWilczynski} An admissible ideal $I\subset 2^{\mathbf{N}}$ is said to satisfy the {\it condition} ({\it AP}) if for every countable family of mutually disjoints sets $\{A_1, A_2, ...\}$ belonging to $I$ there exists a countable family of sets $ \{B_1, B_2,...\}$ such that $A_j \Delta B_j$ is a finite set for $j\in \mathbf{N}$ and $B=\bigcup_{j=1}^{\infty} B_j \in I.$\newline

%\textbf{Definition 1.3.}\cite{KostyrkoSalatWilczynski} A sequence $x=( x_{k})$ of points in $\mathbf{R}$ is said to be $I^{\ast}$-convergent to a real number $\ell$ if there exists a set $M\in F(I)$ (i.e. $\mathbf{N}-M \in I$), $ M=\{k_m: k_1< k_2<...<k_m <...\}$ such that $ \lim_{m}  x_{k_{m}}  =\ell.$ In this case we write $I^{\ast}-\lim x_k = \ell $ and $\ell$ is called the $I^{\ast}-limit$ of $x.$  \newline

{\bf Definition 2.1.} A {\it distribution function} (briefly a d.f.) $F$ is a function from the extended reals
$(-\infty,+\infty)$ into $[0, 1]$ such that
\begin{enumerate}
\item[(a)] it is non-decreasing ;
\item[(b)] it is left-continuous on $(-\infty,+\infty)$;
\item[(c)] $F(-\infty) = 0$ and $F(+\infty) = 1.$
\end{enumerate}

The set of all d.f.'s will be denoted by $\Delta.$ The subset of $\Delta$ consisting of proper d.f's, namely
of those elements $F$ such that $\ell^{+}F(-\infty) = F(-\infty)=0$ and $\ell^{-}F(+\infty) = F(+\infty) =1$ will be
denoted by $D.$ A {\it distance distribution function} (briefly, d.d.f.) is a d.f. $F$ such that $F(0) = 0.$
The set of all d.d.f.f's will be denoted by $\Delta^{+},$ while $D^{+} := D\cap \Delta^{+}$ will denote the set of proper
d.d.f.'s.\\

{\bf Definition 2.2.} A {\it triangular norm} or, briefly, a {\it t-norm} is a binary operation $T : [0, 1]\times [0,1]\rightarrow [0, 1]$ that satisfies the following conditions (see \cite{KlementMesiarPap}):
\begin{enumerate}
\item[(T1)] $T$ is commutative, i.e.,$ T(s, t) = T(t, s)$ for all $s$ and $t$ in $[0, 1];$
\item[(T2)] $T$ is associative, i.e., $T(T(s, t), u) = T(s, T(t, u)) $ for all $s, t$ and $u$ in $[0, 1];$
\item[(T3)] $T$ is nondecreasing, i.e., $T(s, t) \leq T(s', t)$ for all $t \in [0, 1]$ whenever $s \leq s';$
\item[(T4)] $T$ satisfies the boundary condition $T(1, t) = t$ for every $t \in [0, 1].$
\end{enumerate}

$T^{\ast}$ is a continuous {\it t-conorm}, namely, a continuous binary operation on $[0, 1]$ that is related
to a continuous $t$-norm through $T^{\ast}(s, t) = 1- T(1 - s, 1 - t).$ Notice that by virtue of its
commutativity, any $t$-norm $ T$ is nondecreasing in each place. Some examples of $t$-norms $ T$
and its $ t$-conorms $ T^{\ast}$ are: $ M(x, y) = \min\{x, y\}, \Pi(x, y) = x.y$ and $M^{\ast}(x, y) =\max\{x, y\},
\Pi^{\ast}(x, y) = x + y - x . y.$\\

Using the definitions just given above \u{S}herstnev \cite{Sherstnev} defined a PN space as follows:\\

{\bf Definition 2.3.} A triplet $(X, \nu, T)$ is called a {\it probabilistic normed space} (in short PNS) if $X$ is a real vector
space, $ \nu $ is a mapping from $X$ into $D$ and for $x \in X,$ the d.f. $ \nu(x) $ is denoted by $\nu_{x}, \nu_{x}(t) $ is
the value of $\nu_x $ at $t\in \mathbb{ R}$ and $T$ is a {\it t-norm.} $\nu $ satisfies the following conditions :
\begin{enumerate}
\item[(i)] $\nu_{x}(0) = 0;$
\item[(ii)] $\nu_{x}(t) = 1$ for all $t > 0$ if and only if $x = 0;$
\item[(iii)] $\nu_{ax}(t) = \nu_{x}\left(\frac{ t}{|a|}\right )$ for all $ a\in \mathbb{ R} \diagdown \{0\};$
\item[(iv)] $\nu_{x+y}(s + t) \geq  T(\nu_{x}(s), \nu_{y}(t))$ for all $x, y \in X $ and $s, t \in \mathbb{R}^{+}_{0}.$
\end{enumerate}
% = \{x \in \mathbf{R} : x \geq 0\}.$\\

Let $(X, ||.||)$ be a normed space and $\mu\in D$ with $\mu(0)=0$ and $\mu \neq \epsilon_{0},$ where
\[
 \epsilon_{0}(t)=\left\{\begin{array}{cc} 
 0, & \mbox{ if } t\leq 0 \\
 %r-[\sqrt{\lambda_r}]+1\leq k\leq r \mbox{ and } s-[\sqrt{\lambda_s}]+1\leq l\leq s ;\\ 
 1, & \mbox{if } t>0\end{array}\right .\]

For $x\in X, t\in \mathbb{R},$ if we define
\[\nu_{x}(t) = \mu \left(\frac{t}{||x||}\right), x\neq 0,\]
then in \cite{LafuerzaLallenaSempi}, it is proved that $(X, \nu, T)$ is a PN space in the sense of Definition 2.3.
Alsina et al. \cite{AlsinaSchweizerSklar93} gave new definition of a PN-Space. Before giving this, we recall for the
reader's convenience the concept of a triangle function, that of a PN space from the point of
view of the new definition.\\

{\bf Definition 2.4.}  A triangle function is a mapping $\tau$ from $\Delta^{+}\times \Delta^{+}$ into $\Delta^{+} $ such that, for all
$F, G, H, K $ in  $\Delta^{+},$
\begin{enumerate}
\item[(1)] $\tau (F, \epsilon_{0}) = F;$
\item[(2)] $\tau (F,G) = \tau (G, F);$
\item[(3)] $ \tau (F,G) \leq \tau (H,K)$ whenever $F \leq H, ~G \leq K;$
\item[(4)] $\tau (\tau(F,G),H) = \tau (F, \tau(G,H)).$
\end{enumerate}

Particular and relevant triangle functions are the functions $\tau_T , \tau_{T^{\ast}}$ and those of the form
$\Pi_T$ which, for any continuous $t$-norm $T,$ and any $x > 0,$ are given by
\[\tau_{T}(F,G)(x) = \sup\{T(F(u),G(v)) : u + v = x\},\]
\[\tau_{T^{\ast}} (F,G)(x) = \inf\{T^{\ast}(F(u),G(v)) : u + v = x\}\]
and
\[\Pi_{T} (F,G)(x) = T(F(x),G(x)).\]

{\bf Definition 2.5. }(\cite{AlsinaSchweizerSklar93}) A probabilistic normed space is a quad-ruple $(X, \nu, \tau, \tau^{\ast}),$ where $X$ is a
real linear space, $ \tau$ and $\tau^{\ast}$ are continuous triangle functions such that $\tau\leq \tau^{\ast}$ and the mapping
$\nu : X \rightarrow \Delta^{+}$ called the probabilistic norm, satisfies for all $p$ and $q$ in $X,$ the conditions
\begin{enumerate}
\item[(PN1)] $\nu_p = \epsilon_0$ if and only if $ p = \theta$ ($\theta$ is the null vector in $X$);
\item[(PN2)] $\forall  p \in X, \nu_{-p} = \nu_p;$
\item[(PN3)] $\nu_{p+q} \geq \tau(\nu_p, \nu_q);$
\item[(PN4)] $\forall a \in[0, 1] , \nu_p \leq \tau^{\ast}(\nu_{ap}, \nu_{(1-a)p}).$
\end{enumerate}

If a PN space $(X, \nu, \tau, \tau^{\ast}),$ satisfies the following condition

\qquad(\u{S}) $\forall  p \in X, \forall  \lambda\in \mathbb{R}\diagdown\{0\} , \forall t>0 , \nu_{\lambda p}(t) = \nu_{p}\left(\frac{t}{|\lambda|}\right),$
then it is called a \u{S}herstnev PN space; the condition (\u{S}) implies that the best-possible selection
for $\tau^{\ast}$ is $\tau^{\ast} = \tau_{M},$ which satisfies a stricter version of (PN4), namely,

\[\forall a \in[0, 1] , \nu_p = \tau_{M}(\nu_{ap}, \nu_{(1-a)p}).\]

{\bf Definition 2.6.} A Menger PN space under $T$ is a PN space $(X, \nu, \tau, \tau^{\ast})$ denoted by $(X, \nu, T),$
in which $\tau = \tau_{T}$ and $\tau^{\ast} = \tau_{T^{\ast}} ,$ for some continuous $t$-norm $T$ and its $t$-conorm $T^{\ast}.$\\

{\bf Lemma 2.1} (\cite{LafuerzaLallenaSempi}) {\it The simple space generated by $(X, ||.||)$ and by $\mu$ is a Menger PN space
under M and also a \u{S}herstnev PN space. Here $M(x, y) := \min\{x, y\}.$}\\

For further study, by a PN space we mean a PN space in the sense of Definition 2.3. We
now give a quick look on the characterization of convergence and Cauchy sequences on
these spaces. \\

Let $(X, \nu, T)$ be a PN space and $x = (x_k)$ be a sequence in $X.$ We say that $(x_k)$ is
{\it convergent} to $\ell\in X$ with respect to the probabilistic norm $\nu$ if for each $\varepsilon > 0$ and $\alpha\in (0, 1)$
there exists a positive integer $m$ such that $\nu_{x_{k}-\ell}(\varepsilon)> 1-\alpha$ whenever $k\geq m.$ The element
$\ell$ is called the ordinary double limit of the sequence $(x_k)$ and we shall write $\nu-\lim x_k =\ell$ or
$x_k \stackrel{\nu}{\rightarrow} \ell$ as $k\rightarrow \infty.$\\
% with respect to the probabilistic norm ¥í.

A  sequence $(x_k)$ in $X$ is said to be {\it Cauchy} with respect to the probabilistic norm $\nu$
if for each $\varepsilon > 0$ and $\alpha\in (0, 1)$ there exist a positive integer $M = M(\varepsilon,\alpha)$ such that $\nu_{x_{k}- x_{p}}(\varepsilon) > 1 -\alpha$ whenever $ k,p\geq M.$\\

\textbf{Definition 2.7.} Let $\left(  X,\nu,T\right)  $ be an
probabilistic normed space, and let $r\in\left(  0,1\right)  $
and $x\in X.$ The set
\[
B\left(  x,r;t\right)  =\left\{  y\in X:\nu_{y-x}(t) >1-r\right\}
\]
is called open ball with center $x$ and radius $r$ with respect to $t$.\\

Throughout the paper, we denote $I$ is an admissible
ideal of subsets of $%
%TCIMACRO{\U{2115} }%
%BeginExpansion
\mathbb{N}
%EndExpansion
$  and  $\theta=\left(k_{r}\right), $ respectively, unless otherwise stated.

\section{Main results}

 We now obtain our main results.\\

\textbf{Definition 3.1. } Let $I\subset 2^{\mathbb{N}}$ and $\left( X,\nu,T\right)  $ be an PNS. A sequence
$x=\left(  x_{k}\right)  $ in $X$ is said to be $I_{\theta}$-convergent to
$L\in X$ with respect to the probabilistic norm $ \nu$ if, for every $\varepsilon>0 $ and $\alpha \in\left(  0,1\right)  $   the set%
\[
\left\{  r\in%
%TCIMACRO{\U{2115} }%
%BeginExpansion
\mathbb{N}
%EndExpansion
:\text{ }\frac{1}{h_{r}}\sum_{k\in J_{r}}\nu_{  x_{k}-L}(\varepsilon)
\leq 1-\alpha  \right\}  \in I.
\]
$L$ is called the $I_{\theta}-$limit of the sequence $x=\left(  x_{k}\right)  $ in $X,$
and we write $I_{\theta}^{\nu}-\lim x=L.$\\

\textbf{Example 3.1. } Let $\left(
%TCIMACRO{\U{211d} }%
%BeginExpansion
\mathbb{R}
%EndExpansion
,\left\vert .\right\vert \right)  $ denote the space of all real numbers with
the usual norm, and let $T(a,b)=ab$ 
for all $a,b\in\left[  0,1\right]  .$ For all $x\in%
%TCIMACRO{\U{211d} }%
%BeginExpansion
\mathbb{R}
%EndExpansion
$ and every $t>0$, consider $\nu_{x}(t)  =\frac{t}{t+|x|}.$  Then $\left(
%TCIMACRO{\U{211d} }%
%BeginExpansion
\mathbb{R},
%EndExpansion
\nu,T\right)  $ is an PNS. If we take $I=\left\{  A\subset%
%TCIMACRO{\U{2115} }%
%BeginExpansion
\mathbb{N}
%EndExpansion
:\text{ }\delta\left(  A\right)  =0\right\}  ,$ where $\delta\left(  A\right)
$ denotes the natural density of the set A, then $I$ is a non-trivial
admissible ideal. Define a sequence $x=\left(  x_{k}\right)  $ as follows:%
 \[x_{k}=\left\{\begin{array}{cc} 
 1, & \mbox{ if } k=i^{2}, i\in\mathbb{N}\\
 %r-[\sqrt{\lambda_r}]+1\leq k\leq r \mbox{ and } s-[\sqrt{\lambda_s}]+1\leq l\leq s ;\\ 
 0, & \mbox{otherwise. }\end{array}\right .\]

%\[x_{k}=\{%
%TCIMACRO{\QATOP{1,\text{ if }k=i^{2}\text{ }\left(  i\in\U{2115} \right)
%}{0,\text{ otherwise \ \ \ \ \ \ \ \ }}}%
%BeginExpansion
%\genfrac{}{}{0pt}{}{(1,0),\text{ if }k=i^{2}\text{ }\left(  i\in\mathbb{N} \right)  }{(0,),\text{ otherwise \ \ \ \ \ \ \ \ }}%
%EndExpansion .\]
Then for every $\alpha\in\left(  0,1\right)  $ and for any $\varepsilon>0,$ the set%
\[
K  =\left\{  r\in%
%TCIMACRO{\U{2115} }%
%BeginExpansion
\mathbb{N}
%EndExpansion
:\text{ }\frac{1}{h_{r}}\sum_{k\in J_{r}}\nu_{ x_{k}}(\varepsilon)
\leq1-\alpha\right\}
\]
will be a finite set. Hence, $\delta\left( K 
\right)  =0$ and  consequently $K \in I,$ i.e.,
$I_{\theta}^{\nu}-\lim x=0.$\\

\textbf{Lemma 3.1.} {\it Let $\left(  X,\nu,T\right)  $ be an PNS
and $x=\left(  x_{k}\right)  $ be a sequence in $X.$ Then, for every
$\varepsilon>0$ and $\alpha\in (0,1)$ the following statements are equivalent:}

\qquad$(i)$ $I_{\theta}^{\nu}-\lim x=L,$

\qquad$\left(  ii\right)  $ $\left\{  r\in%
%TCIMACRO{\U{2115} }%
%BeginExpansion
\mathbb{N}
%EndExpansion
:\text{ }\frac{1}{h_{r}}\sum_{k\in J_{r}}\nu_{  x_{k}-L} (\varepsilon)\leq1-\alpha\text{ }\right\}  \in I$

%and $\left\{  n\in%
%TCIMACRO{\U{2115} }%
%BeginExpansion
%\mathbb{N}
%EndExpansion
%:\text{ }\frac{1}{\lambda_{n}}\sum_{k\in J_{n}}\nu\left(  x_{k}-L,y;t\right)
%\geq\varepsilon\text{ }\right\}  \in I,$

\qquad$\left(  iii\right)  $ $\left\{
\begin{array}
[c]{c}r\in%
%TCIMACRO{\U{2115} }%
%BeginExpansion 
\mathbb{N}
%EndExpansion 
:\text{ }\frac{1}{h_{r}}\sum_{k\in J_{r}}\nu_{  x_{k}-L} (\varepsilon)>1-\alpha
%\text{ }\\ \text{and }\frac{1}{\lambda_{n}}\sum_{k\in J_{n}}\nu\left(  x_{k} -L,y;t\right)  <\varepsilon\text{ } 
\end{array} \right\}  \in F\left(  I\right)  ,$

%\qquad$\left(  iv\right)  $ $\left\{  n\in%
%TCIMACRO{\U{2115} }%
%BeginExpansion
%\mathbb{N}
%EndExpansion
%:\text{ }\frac{1}{\lambda_{n}}\sum_{k\in J_{n}}\nu_{ x_{k}-L}(\varepsilon)
%>1-\alpha\text{ }\right\}  \in F\left(  I\right)  $

%and $\left\{  n\in%
%TCIMACRO{\U{2115} }%
%BeginExpansion
%\mathbb{N}
%EndExpansion
%:\text{ }\frac{1}{\lambda_{n}}\sum_{k\in J_{n}}\nu\left(  x_{k}-L,y;t\right)
%<\varepsilon\text{ }\right\}  \in F\left(  I\right)  ,$

\qquad$\left(  iii\right)  $ $I_{\theta}-\lim \nu_{  x_{k}-L}(\varepsilon)  =1.$\\
%and $I_{\lambda}-\lim\nu\left(  x_{k}-L,y;t\right)  =0.$\\

\textbf{Theorem 3.1.} {\it Let $\left(  X,\nu,T\right)  $ be an
PNS and if a sequence $x=\left(  x_{k}\right)  $ in $X$ is $I_{\theta}%
$-convergent to $L\in X$ with respect to the probabilistic norm
$\nu,$ then $I_{\theta}^{\nu}-\lim x$ is unique.}\\

\textbf{Proof.} Suppose that $I_{\theta}^{\nu}-\lim x=L_{1}$ and $I_{\theta}^{\nu}-\lim x=L_{2}$
$\left(  L_{1}\neq L_{2}\right).$ Given $\alpha>0 $ and choose $\beta\in\left(  0,1\right)  $
%and  $\beta\in\left(  0,1\right)  $ 
such that 
\begin{equation}
T(  1-\beta, 1-\beta)
  >1-\alpha.
  \end{equation} 
Then for $\varepsilon>0$, define the following sets:%
\[
K_{1}  =\left\{  r\in%
%TCIMACRO{\U{2115} }%
%BeginExpansion
\mathbb{N}
%EndExpansion
:\text{ }\frac{1}{h_{r}}\sum_{k\in J_{r}}\nu_{x_{k}-L_{1}}\left(\frac{\varepsilon}{2}\right)  \leq1-\beta\text{ }\right\}  ,
\]

\[
K_{2} =\left\{  r\in%
%TCIMACRO{\U{2115} }%
%BeginExpansion
\mathbb{N}
%EndExpansion
:\text{ }\frac{1}{h_{r}}\sum_{k\in J_{r}}\nu_{x_{k}-L_{2}}\left(\frac{\varepsilon}{2}\right)  \leq1-\beta\text{ }\right\}  ,
\]%
%\[K_{\nu,1}\left(  \beta,t\right)  =\left\{  n\in%
%TCIMACRO{\U{2115} }%
%BeginExpansion
%\mathbb{N}
%EndExpansion
%:\text{ }\frac{1}{\lambda_{n}}\sum_{k\in J_{n}}\nu\left(  x_{k}-L_{1}%
%,y;\frac{t}{2}\right)  \geq\beta\text{ }\right\}  ,
%\]%
%\[K_{\nu,2}\left(  \beta,t\right)  =\left\{  n\in%
%TCIMACRO{\U{2115} }%
%BeginExpansion
%\mathbb{N}
%EndExpansion
%:\text{ }\frac{1}{\lambda_{n}}\sum_{k\in J_{n}}\nu\left(  x_{k}-L_{2}%
%,y;\frac{t}{2}\right)  \geq\beta\text{ }\right\}  .\]
Since $I_{\theta}^{\nu}-\lim x=L_{1}$, using Lemma 2.1.,
we have $K_{1} \in I.$  Also, using $I_{\theta}^{\nu }-\lim x=L_{2},$ we get $K_{2} \in
I.$  Now let%
\[
K = K_{1}
\cup K_{2} .
\]
Then $K \in I.$ This implies that its
complement $K^{c} $ is a non-empty set in
$F(I)$. Now if $r\in K^{c} ,$  let us
consider  $r\in  K_{1}^{c}  \cap
K_{2}^{c} .$ Then we have%
\[
\text{ }\frac{1}{h_{r}}\sum_{k\in J_{r}}\nu_{ x_{k}-L_{1}}\left(\frac{\varepsilon}{2}\right)  >1-\beta\text{ and }\frac{1}{h_{r}}\sum_{k\in
J_{r}}\nu_{x_{k}-L_{2}}\left(\frac{\varepsilon}{2}\right)  >1-\beta\text{.}%
\]
Now, we choose a $s\in%
%TCIMACRO{\U{2115} }%
%BeginExpansion
\mathbb{N}
%EndExpansion
$ such that%
\[
\nu_{ x_{s}-L_{1}}\left(\frac{\varepsilon}{2}\right)  >\frac{1}{h_{r}}\sum_{k\in
J_{r}}\nu_{ x_{k}-L_{1}}\left(\frac{\varepsilon}{2}\right)  >1-\beta\text{ }%
\]
and%
\[
\nu_{  x_{s}-L_{2}}\left(\frac{\varepsilon}{2}\right)  >\frac{1}{h_{r}}\sum_{k\in
J_{r}}\nu_{ x_{k}-L_{2}}\left(\frac{\varepsilon}{2}\right)  >1-\beta\text{ }%
\]
e.g., consider $\max\left\{  \nu_{ x_{k}-L_{1}}\left(\frac{\varepsilon}{2}\right)
,\nu_{ x_{k}-L_{2}}\left(\frac{\varepsilon}{2}\right)  :\text{ }k\in J_{r}\right\}  $
and choose that $k$ as $s$ for which the maximum occurs. Then from (2.1), we have%
\[
\nu_{ L_{1}-L_{2}}(\varepsilon)  \geq T\left(\nu_{ x_{s}-L_{1}}\left(\frac{\varepsilon}%
{2}\right) ,\nu_{x_{s}-L_{2}}\left(\frac{\varepsilon}{2}\right)\right)
>T\left(  1-\beta,  1-\beta \right) >1-\alpha.
\]
Since $\alpha>0$ is arbitrary, we have $\nu_{  L_{1}-L_{2}}(\varepsilon)  =1$ for all $\varepsilon>0$, which implies that $L_{1}=L_{2}.$ 
%On the other hand, if $n\in\left[  K_{\nu,1}^{C}\left(  \beta,t\right)  \cap K_{\nu,2} ^{C}\left(  \beta,t\right)  \right]  ,$ then using a similar technique, it can be proved that $\nu\left(  L_{1}-L_{2},y;t\right)  <\varepsilon$ for all $t>0$ and arbitrary $\varepsilon>0.$ Thus we obtain $L_{1}=L_{2}.$ 
Therefore, we
conclude that $I_{\theta}^{\nu}-\lim x$ is unique.\\

 Here, we introduce the notion of $\theta$-convergence in an PNS and
discuss some properties.\\

\textbf{Definition 3.2.} Let $\left(  X,\nu,T\right)  $ be an
PNS. A sequence $x=\left(  x_{k}\right)  $ in $X$ is $\theta$-convergent
to $L\in X$ with respect to the probabilistic norm $\nu$ if, for $\alpha\in\left(  0,1\right)  $ and every
$\varepsilon>0$, there exists $r_{o}\in%
%TCIMACRO{\U{2115} }%
%BeginExpansion
\mathbb{N}
%EndExpansion
$ such that%
\[
\frac{1}{h_{r}}\sum_{k\in J_{r}}\nu_{  x_{k}-L}(\varepsilon)
>1-\alpha
\]
for all $r\geq r_{o}.$ In this case, we write $\nu^{\theta}-\lim x=L.$\\

\textbf{Theorem 3.2.} {\it Let $\left(  X,\nu,T\right)  $ be an
PNS and let $x=\left(  x_{k}\right)  $ in $X$. If $x=\left(  x_{k}\right)$ is $\theta$-convergent with respect to the probabilistic norm
$\nu,$ then $\nu^{\theta}-\lim x$ is unique.}\\

\textbf{Proof. } Suppose that $\nu^{\theta}-\lim x=L_{1}$ and $\nu^{\theta}-\lim x=L_{2}$
$\left(  L_{1}\neq L_{2}\right)  .$ Given $\alpha\in\left(  0,1\right)  $
and choose $\beta\in\left(  0,1\right)  $ such that $T\left(  1-\beta,  1-\beta\right)  >1-\alpha.$ 
%and $\beta o\beta<\varepsilon.$
Then for any $\varepsilon>0$, there exists $r_{1}\in%
%TCIMACRO{\U{2115} }%
%BeginExpansion
\mathbb{N}
%EndExpansion
$ such that%
\[
\frac{1}{h_{r}}\sum_{k\in J_{r}}\nu_{ x_{k}-L_{1}}\left(\varepsilon\right)
>1-\alpha
\]
for all $r\geq r_{1}.$ Also, there exists $r_{2}\in%
%TCIMACRO{\U{2115} }%
%BeginExpansion
\mathbb{N}
%EndExpansion
$ such that%
\[
\frac{1}{h_{r}}\sum_{k\in J_{r}}\nu_{ x_{k}-L_{2}}\left(\varepsilon\right)
>1-\alpha
\]
for all $r\geq r_{2}.$ Now, consider $r_{o}=\max\left\{  r_{1},r_{2}\right\}
.$ Then for $r\geq r_{o}$, we will get a $s\in%
%TCIMACRO{\U{2115} }%
%BeginExpansion
\mathbb{N}
%EndExpansion
$ such that%
\[
\nu_{ x_{s}-L_{1}}\left(\frac{\varepsilon}{2}\right)  >\frac{1}{h_{r}}\sum_{k\in
J_{r}}\nu_{ x_{k}-L_{1}}\left(\frac{\varepsilon}{2}\right)  >1-\beta\text{ }%
\]
and%
\[
\nu_{  x_{s}-L_{2}}\left(\frac{\varepsilon}{2}\right)  >\frac{1}{h_{r}}\sum_{k\in
J_{r}}\nu_{ x_{k}-L_{2}}\left(\frac{\varepsilon}{2}\right)  >1-\beta.\text{ }%
\]
Then, we have%
\[
\nu_{L_{1}-L_{2}}(\varepsilon)  \geq T\left(\nu_{ x_{s}-L_{1}}\left(\frac{\varepsilon}%
{2}\right)  , \nu_{ x_{s}-L_{2}}\left(\frac{\varepsilon}{2}\right)\right)
>T(  1-\beta,  1-\beta) >1-\alpha.
\]
Since $\alpha>0$ is arbitrary, we have $\nu_{ L_{1}-L_{2}%
}(\varepsilon) =1$ for all $\varepsilon>0,$
%, by using a similar technique, it can be proved that $\nu\left(  L_{1}-L_{2},y;t\right)  =0$ for all $t>0,$
 which implies that $L_{1}=L_{2}.$\\
 
 \textbf{Theorem 3.3.} {\it Let $\left(  X,\nu,T\right)  $ be an
PNS and let $x=\left(  x_{k}\right)  $ in $X$. If $\nu^{\theta}-\lim x=L,$ then $I_{\theta}^{\nu}-\lim x=L.$}\\

\textbf{Proof.} Let $\nu^{\theta}-\lim x=L,$
then for every $\varepsilon>0$ and given $\alpha\in\left(  0,1\right)  $, there
exists $r_{0}\in%
%TCIMACRO{\U{2115} }%
%BeginExpansion
\mathbb{N}
%EndExpansion
$ such that%
\[
\frac{1}{h_{r}}\sum_{k\in J_{r}}\nu_{ x_{k}-L}(\varepsilon)
>1-\alpha
\]
for all $r\geq r_{0}.$ Therefore the set%
\[
B=\left\{  r\in%
%TCIMACRO{\U{2115} }%
%BeginExpansion
\mathbb{N}
%EndExpansion
:\text{ }\frac{1}{h_{r}}\sum_{k\in J_{r}}\nu_{ x_{k}-L}(\varepsilon)
\leq 1-\alpha  \right\}
\subseteq\left\{  1,2,...,n_{0}-1\right\}  .
\]
But, with $I$ being admissible, we have $B\in I.$ Hence $I_{\theta}^{\nu}-\lim x=L.$\\

\textbf{Theorem 3.4.} {\it Let $\left(  X,\nu,T\right)  $ be an
PNS and   $x=\left(  x_{k}\right) , y=(y_k) $ be two sequence in $X.$}\\

(i) {\it If $I_{\theta}^{\nu}-\lim x_k =L_1$ and $I_{\theta}^{\nu}-\lim y_k =L_2,$ then $I_{\theta}^{\nu}-\lim(x_k \pm y_k) =L_1 \pm L_2;$}\\

(ii) {\it If $I_{\theta}^{\nu}-\lim x_k =L$ and $ a$ be a non-zero real number, then $I_{\theta}^{\nu}-\lim a x_k=aL.$ If $a=0,$ then  result is true only if $I$ is an  admissble of $\mathbb{N}$.}\\

{\bf Proof.} (i) We have proved that, if $I_{\theta}^{\nu}-\lim x_k =L_1$ and $I_{\theta}^{\nu}-\lim y_k =L_2,$ then $I_{\theta}^{\nu}-\lim(x_k + y_k) =L_1 + L_2,$ only. The proof of the other partfollows similarly.\\

Take $\varepsilon>0, \alpha\in (0,1) $ and choose $\beta\in\left(  0,1\right) $ such that the condition (3.1) holds. If we define
\[
A_{1}  =\left\{  r\in \mathbb{N}
:\text{ }\frac{1}{h_{r}}\sum_{k\in J_{r}}\nu_{x_{k}-L_{1}}\left(\frac{\varepsilon}{2}\right)  \leq1-\beta\text{ }\right\}  ,
\]

and
\[
A_{2}  =\left\{  r\in \mathbb{N}
:\text{ }\frac{1}{h_{r}}\sum_{k\in J_{r}}\nu_{y_{k}-L_{2}}\left(\frac{\varepsilon}{2}\right)  \leq1-\beta\text{ }\right\} .
\]
Then $A_1^{c}\cap A_2^{c}\in F(I).$ We claim that 
\[A_1^{c}\cap A_2^{c}\subset \left\{  r\in \mathbb{N}
:\text{ }\frac{1}{h_{r}}\sum_{k\in J_{r}}\nu_{(x_k - L_1)+( y_{k}-L_{2})}(\varepsilon)  >1-\alpha\right\}.\]
 Let $n\in A_1^{c}\cap A_2^{c}.$ Now, using (3.1), we have
 \[\frac{1}{h_{r}}\sum_{n\in J_{r}}\nu_{(x_n -L_1)+ (y_n - L_2)}(\varepsilon)\geq T \left(\frac{1}{h_{r}}\sum_{n\in J_{r}}\nu_{x_{n}-L_{1}}\left(\frac{\varepsilon}{2}\right),\frac{1}{h_{r}}\sum_{n\in J_{r}}\nu_{y_{n}-L_{2}}\left(\frac{\varepsilon}{2}\right)\right) \]
 \[> T(1-\beta, 1-\beta)> 1-\alpha.\]
 Hence
 \[A_1^{c}\cap A_2^{c} \subset \left\{  r\in \mathbb{N}
:\text{ }\frac{1}{h_{r}}\sum_{k\in J_{r}}\nu_{(x_k - L_1)+( y_{k}-L_{2})}(\varepsilon)  >1-\alpha\right\}.\]
As $A_1^{c}\cap A_2^{c}\in F(I),$ so 
\[\left\{  r\in \mathbb{N}
:\text{ }\frac{1}{h_{r}}\sum_{k\in J_{r}}\nu_{(x_k - L_1)+( y_{k}-L_{2})}(\varepsilon)  \leq 1-\alpha\right\}\in I.\]
Therefore $I_{\theta}^{\nu}-\lim(x_k + y_k) =L_1 + L_2.$\\

(ii) Suppose $a\neq 0.$ Since $I_{\theta}^{\nu}-\lim x_k =L,$ for each $\varepsilon>0$ and $ \alpha\in (0,1), $ the set 
\[A(\varepsilon, \alpha)=\left\{r\in \mathbb{N}: \frac{1}{h_r}\sum_{k\in J_{r}}\nu_{x_k - L}(\varepsilon)  < 1-\alpha\right\}\in F(I).\]
If $n\in A(\varepsilon, \alpha),$  the we have
\[\frac{1}{h_r}\sum_{k\in J_{r}}\nu_{ax_k -a L}(\varepsilon)= \frac{1}{h_r}\sum_{k\in J_{r}}\nu_{x_k - L}\left(\frac{\varepsilon}{|a|}\right)\]
\[\geq T \left(\frac{1}{h_r}\sum_{k\in J_{r}}\nu_{x_k - L}(\varepsilon), \nu_{0}\left(\frac{\varepsilon}{|a|}-\varepsilon\right)\right)\]
\[\geq T\left(\frac{1}{h_r}\sum_{k\in J_{r}}\nu_{x_k - L}(\varepsilon),1\right)
\geq \frac{1}{h_r}\sum_{k\in J_{r}}\nu_{x_k - L}(\varepsilon)> 1-\alpha\]
Hence 
\[A(\varepsilon, \alpha)\subset \left\{r\in \mathbb{N}: \frac{1}{h_r}\sum_{k\in J_{r}}\nu_{ax_k -a L}(\varepsilon)> 1-\alpha\right\}\]
and
\[\left\{r\in \mathbb{N}: \frac{1}{h_r}\sum_{k\in J_{r}}\nu_{ax_k -a L}(\varepsilon)> 1-\alpha\right\}\in F(I).\]
It follows that
\[\left\{r\in \mathbb{N}: \frac{1}{h_r}\sum_{k\in J_{r}}\nu_{ax_k -a L}(\varepsilon)\leq 1-\alpha\right\}\in I.\]
Hence $I_{\theta}^{\nu}-\lim a x_k=aL.$\\

Next suppose that $a=0.$ Then for each $\varepsilon>0$ and $ \alpha\in (0,1), $ we have
\[\frac{1}{h_r}\sum_{k\in J_{r}}\nu_{0x_k -0 L}(\varepsilon)=\frac{1}{h_r}\sum_{k\in J_{r}}\nu_{0}(\varepsilon)=1 > 1-\alpha,\]
it follows that $\nu^{\theta}-\lim x=\ell.$ Hence from Theorem 3.3, $I_{\theta}^{\nu}-\lim x = \ell.$\\

%\textbf{Theorem 3.4.} {\it Sequential method $I_{\theta}$ is regular. i.e. If $\nu^{\theta}-\lim x=\ell,$ then $I_{\theta}^{\nu}-\lim x = \ell.$}\\

%\textbf{Proof.} The proof follows from the fact that $I$ is admissible and Theorem 3.3 (A).\\

\textbf{Theorem 3.5. } {\it Let $\left(  X,\nu,T\right)  $ be an
PNS and let $x=\left(  x_{k}\right)  $ in $X$. If $\nu^{\theta}-\lim x=L,$ then there exists a subsequence $\left(  x_{m_{k}%
}\right)  $ of $x=\left(  x_{k}\right)  $ such that $\nu-\lim x_{m_{k}}=L.$}\\

\textbf{Proof.} Let $\nu^{\theta}-\lim x=L.$
Then, for every $\varepsilon>0$ and given $\alpha\in\left(  0,1\right)  $, there
exists $r_{0}\in%
%TCIMACRO{\U{2115} }%
%BeginExpansion
\mathbb{N}
%EndExpansion
$ such that%
\[
\frac{1}{h_{r}}\sum_{k\in J_{r}}\nu_{ x_{k}-L}(\varepsilon)
>1-\alpha
\]
for all $r\geq r_{0}.$ Clearly, for each $r\geq r_{0},$ we can select an
$m_{k}\in J_{r}$ such that%
\[
\nu_{x_{m_{k}}-L}(\varepsilon)  >\frac{1}{h_{r}}\sum_{k\in J_{r}}%
\nu_{ x_{k}-L}(\varepsilon)  >1-\alpha.
\]
It follows that $\nu-\lim x_{m_{k}}=L.$\\

\textbf{Definition 3.3.} Let $\left(  X,\nu,T\right)  $ be an
PNS and let $x=\left(  x_{k}\right)  $ be a sequence in $X$. Then,
\begin{enumerate}
\item[(1)] An element $L\in X$ is said to be $I_{\theta}$-limit point of
$x=\left(  x_{k}\right)  $ if there is a set $M=\left\{  m_{1}<m_{2}%
<...<m_{k}<...\right\}  \subset%
%TCIMACRO{\U{2115} }%
%BeginExpansion
\mathbb{N}
%EndExpansion
$ such that the set $M^{\imath}=\left\{  r\in%
%TCIMACRO{\U{2115} }%
%BeginExpansion
\mathbb{N}
%EndExpansion
:\text{ }m_{k}\in J_{r}\right\}  \notin I$ and $\nu^{\theta}-\lim x_{m_{k}}=L.$

\item[(2)] An element $L\in X$ is said to be $I_{\theta}$-cluster point of
$x=\left(  x_{k}\right)  $ if for every $\varepsilon>0$ and $\alpha\in\left(
0,1\right)  ,$ we have%
\[
\left\{  r\in%
%TCIMACRO{\U{2115} }%
%BeginExpansion
\mathbb{N}
%EndExpansion
:\text{ }\frac{1}{h_{r}}\sum_{k\in J_{r}}\nu_{ x_{k}-L}(\varepsilon)>1-\alpha\right\}  \notin I.
\]
\end{enumerate}

Let $\Lambda_{\nu}^{I_{\theta}}\left(  x\right)$ denote the set of all $I_{\theta}$-limit points and $\Gamma_{\nu}^{I_{\theta}}\left(  x\right)  $ denote the set of all
$I_{\theta}$-cluster points in $X$, respectively.\\

\textbf{Theorem 3.6.} {\it Let $\left(  X,\nu,T\right)  $ be an
PNS. For each sequence $x=\left(  x_{k}\right)  $ in $X$, we have
$\Lambda_{\nu}^{I_{\theta}}\left(  x\right)
\subset\Gamma_{\nu}^{I_{\theta}}\left(  x\right)  .$}\\

%\qquad
\textbf{Proof.} Let $L\in\Lambda_{\nu}^{I_{\theta}}\left(  x\right)  ,$ then there exists a set $M\subset%
%TCIMACRO{\U{2115} }%
%BeginExpansion
\mathbb{N}
%EndExpansion
$ such that $M^{\imath}\notin I,$ where $M$ and $M^{\imath}$ are as in the
Definition 2.3., satisfies $\nu^{\lambda}-\lim
x_{m_{k}}=L.$ Thus, for every $\varepsilon>0$ and $\alpha\in\left(  0,1\right)  ,$
there exists $r_{0}$ $\in%
%TCIMACRO{\U{2115} }%
%BeginExpansion
\mathbb{N}
%EndExpansion
$ such that%
\[
\frac{1}{h_{r}}\sum_{k\in J_{r}}\nu_{ x_{m_{k}}-L}(\varepsilon)
>1-\alpha
\]
for all $r\geq r_{0}.$ Therefore,%
%\begin{align*}
\[
B   =\left\{  r\in%
%TCIMACRO{\U{2115} }%
%BeginExpansion
\mathbb{N}
%EndExpansion
:\text{ }\frac{1}{h_{r}}\sum_{k\in J_{r}}\nu_{x_{k}-L}(\varepsilon)
>1-\alpha\right\} 
 \supseteq M^{\imath}\setminus\left\{  m_{1},m_{2},...,m_{n_{0}}\right\}  .\]
%\end{align*}
Now, with $I$ being admissible, we must have $M^{\imath}\setminus\left\{
m_{1},m_{2},...,m_{k_{0}}\right\}  \notin I$ and as such $B\notin I.$ Hence
$L\in\Gamma_{\nu}^{I_{\theta}}\left(  x\right)  .$\\

\textbf{Theorem 3.7.} {\it Let $\left(  X,\nu,T\right)  $ be an
PNS. For each sequence $x=\left(  x_{k}\right)  $ in $X$, the set
$\Gamma_{\nu}^{I_{\theta}}\left(  x\right)  $ is
closed set in $X$ with respect to the usual topology induced by the
probabilistic norm $\nu^{\theta}.$}\\

\textbf{Proof.} Let $y\in\overline{\Gamma_{\nu}^{I_{\theta}}\left(  x\right)  }.$ Take $\varepsilon>0$ and $\alpha\in\left(
0,1\right)  .$ Then there exists $L_{0}\in \Gamma_{\nu}^{I_{\theta}}\left(  x\right)  \cap B\left(  y,\alpha,\varepsilon\right)  .$
Choose $\delta>0$ such that $B\left(  L_{0},\delta,\varepsilon\right)  \subset B\left(
y,\alpha,\varepsilon\right)  .$ We have%
\[
G=\left\{  r\in%
%TCIMACRO{\U{2115} }%
%BeginExpansion
\mathbb{N}
%EndExpansion
:\text{ }\frac{1}{h_{r}}\sum_{k\in J_{r}}\nu_{ x_{k}-y}(\varepsilon)
>1-\alpha\right\}
\]%
\[
\supseteq\left\{  r\in%
%TCIMACRO{\U{2115} }%
%BeginExpansion
\mathbb{N}
%EndExpansion
:\text{ }\frac{1}{h_{r}}\sum_{k\in J_{r}}\nu_{ x_{k}-L_{0}}(\varepsilon)%
  >1-\delta\right\}=H.
\]
Thus $H\notin I$ and so $G\notin I.$ Hence $y\in \Gamma_{\nu}^{I_{\theta}}\left(  x\right)  .$\\

\qquad\textbf{Theorem 3.8.} {\it Let $\left(  X,\nu,T\right)  $ be an
PNS and let $x=\left(  x_{k}\right)  $ in $X$. Then the following
statements are equivalent:}

\qquad$\left(  1\right)  $ $L$ {\it is a $I_{\theta}-$limit point of $x,$}

\qquad$\left(  2\right)  $ {\it There exist two sequences $y$ and $z$ in $X$ such
that $x=y+z$ and $\nu^{\theta}-\lim y=L$ and
$\left\{  r\in%
%TCIMACRO{\U{2115} }%
%BeginExpansion
\mathbb{N}
%EndExpansion
:\text{ }k\in J_{r},\text{ }z_{k}\neq\overline{\theta}\right\}  \in I,$ where $\overline{\theta}$ is
the zero element of $X$.}\\

\textbf{Proof.} Suppose that (1) holds. Then there exist sets $M$ and
$M^{\imath}$ as in Definition 2.3. such that $M^{\imath}\notin I$ and $\nu^{\theta}-\lim x_{m_{k}}=L.$ Define the sequences $y$ and
$z$ as follows:%

 \[
 y_{k}=\left\{\begin{array}{cc} 
 x_k, & \mbox{ if } k\in J_{r}\text{; }r\in M^{\imath}, \\
 %r-[\sqrt{\lambda_r}]+1\leq k\leq r \mbox{ and } s-[\sqrt{\lambda_s}]+1\leq l\leq s ;\\ 
 L, & \mbox{ otherwise.} \end{array}\right .\]

%\[
%y_{k}=\{%
%TCIMACRO{\QATOP{x_{k}\text{ \ if }k\in J_{n}\text{; }n\in M^{\imath}}{L\text{,
%\ otherwise \ \ \ \ \ \ \ \ \ \ \ }}}%
%BeginExpansion
%\genfrac{}{}{0pt}{}{x_{k}\text{ \ if }k\in J_{n}\text{; }n\in M^{\imath
%}}{L\text{, \ otherwise \ \ \ \ \ \ \ \ \ \ \ }}%
%EndExpansion\]
and%

\[
 z_{k}=\left\{\begin{array}{cc} 
 \overline{\theta}, & \mbox{ if } k\in J_{r}\text{; } r\in M^{\imath}, \\
 %r-[\sqrt{\lambda_r}]+1\leq k\leq r \mbox{ and } s-[\sqrt{\lambda_s}]+1\leq l\leq s ;\\ 
 x_{k}-L, & \mbox{ otherwise.} \end{array}\right .\]

%\[z_{k}=\{%
%TCIMACRO{\QATOP{\theta\text{ \ if }k\in J_{n}\text{; }n\in M^{\imath}}%
%{x_{k}-L\text{, \ otherwise \ \ \ \ \ \ \ \ \ \ \ }}}%
%BeginExpansion 
%\genfrac{}{}{0pt}{}{\theta\text{ \ if }k\in J_{n}\text{; }n\in M^{\imath
%}}{x_{k}-L\text{, \ otherwise \ \ \ \ \ \ \ \ \ \ \ }}%
%EndExpansion.
%\]
It sufficies to consider the case $k\in J_{r}$ such that\ $r\in%
%TCIMACRO{\U{2115} }%
%BeginExpansion
\mathbb{N}
%EndExpansion
\diagdown M^{\imath}.$ Then for each $\alpha\in\left(  0,1\right)  $ and
$\varepsilon>0$, we have $\nu_{ y_{k}-L}(\varepsilon)  =1>1-\alpha.$  Thus, in this case,%
\[
\text{ }\frac{1}{h_{r}}\sum_{k\in J_{r}}\nu_{ y_{k}-L}(\varepsilon)
=1>1-\alpha .
\]
Hence $\nu^{\theta}-\lim y=L.$ Now $\left\{  r\in%
%TCIMACRO{\U{2115} }%
%BeginExpansion
\mathbb{N}
%EndExpansion
:\text{ }k\in J_{r},\text{ }z_{k}\neq\theta\right\}  \subset%
%TCIMACRO{\U{2115} }%
%BeginExpansion
\mathbb{N}
%EndExpansion
\diagdown M^{\imath}$ and so $\left\{  r\in%
%TCIMACRO{\U{2115} }%
%BeginExpansion
\mathbb{N}
%EndExpansion
:\text{ }k\in J_{r},\text{ }z_{k}\neq\theta\right\}  \in I.$

 Now, suppose that (2) holds. Let $M^{\imath}=\left\{  r\in%
%TCIMACRO{\U{2115} }%
%BeginExpansion
\mathbb{N}
%EndExpansion
:\text{ }k\in J_{r},\text{ }z_{k}=\theta\right\}  .$ Then, clearly $M^{\imath
}\in F\left(  I\right)  $ and so it is an infinite set. Construct the set
$M=\left\{  m_{1}<m_{2}<...<m_{k}<...\right\}  \subset%
%TCIMACRO{\U{2115} }%
%BeginExpansion
\mathbb{N}
%EndExpansion
$ such that $m_{k}\in J_{r}$ and $z_{m_{k}}=\overline{\theta}.$ Since $x_{m_{k}}%
=y_{m_{k}}$ and $\nu^{\theta}-\lim y=L$ we obtain
$\nu^{\theta}-\lim x_{m_{k}}=L.$ This completes the proof.\\

\textbf{Theorem 3.9.} \textit{Let $\left(  X,\nu,T\right)  $ be an PNS and $x=(x_{k})$ be a sequence in X. Let $I$ be an admissible ideal in $\mathbb{N}.$ If there is a $I_{\theta}^{\nu}$-convergent sequence $y=(y_{k})$ in X such that } $\{k\in \mathbb{N}: y_{k} \neq x_{k}\}\in I$ {\it then $x$ is also $I_{\theta}^{\nu}$-convergent.}\\

\textbf{Proof.} Suppose that $\{k\in \mathbb{N}: y_{k} \neq x_{k}\}\in I$ and $I_{\theta}^{\nu}-\lim y = \ell.$  
Then for every $\alpha\in\left(  0,1\right)  $ and $\varepsilon>0,$ the set%
\[
\left\{   r\in%
%TCIMACRO{\U{2115} }%
%BeginExpansion
\mathbb{N}
%EndExpansion
:\text{ }\frac{1}{h_{r}}\sum_{k\in J_{r}}\nu_{ y_{k}-L}(\varepsilon)  \leq1-\alpha\right\}  \in
I.
\]

%Then for $t>0$ and for non zero $z\in X$ we get 
%\[ \left\{k\in \mathbf{N}: \mathcal{F}(\Delta^{n}y_k-\ell, z;\frac{t}{2})\leq 1-\varepsilon\right \}\in I.\]

For every $0< \alpha<1$  and $\varepsilon>0,$ we have 

\begin{equation}
\left\{   r\in%
%TCIMACRO{\U{2115} }%
%BeginExpansion
\mathbb{N}
%EndExpansion
:\text{ }\frac{1}{h_{r}}\sum_{k\in J_{r}}\nu_{ x_{k}-L}(\varepsilon)  \leq1-\alpha\right\}
\end{equation}
\[ \subseteq  \{k\in \mathbb{N}: y_{k} \neq x_{k}  \} \cup 
\left\{   r\in%
%TCIMACRO{\U{2115} }%
%BeginExpansion
\mathbb{N}
%EndExpansion
:\text{ }\frac{1}{h_{r}}\sum_{k\in J_{r}}\nu_{y_{k}-L}(\varepsilon)  \leq 1-\alpha\right\} .\]

As both the sets of right-hand side of (2.2)  is in $I,$ therefore we have that
\[\left\{   r\in%
%TCIMACRO{\U{2115} }%
%BeginExpansion
\mathbb{N}
%EndExpansion
:\text{ }\frac{1}{h_{r}}\sum_{k\in J_{r}}\mu\left(  x_{k}%
-L,t\right)  \leq1-\varepsilon\text{ or }\frac{1}{h_{r}}\sum_{k\in
J_{r}}\nu\left(  x_{k}-L,t\right)  \geq\varepsilon\text{ }\right\}\in I.\]  
This completes the proof of the theorem.\\

\textbf{Definition 3.4.} Let $\left(  X,\nu,T\right)  $ be an
PNS. A sequence $x=\left(  x_{k}\right)  $ in $X$ is said to be $\theta$-Cauchy
sequence with respect to the probabilistic norm $\nu$ if, for every $\varepsilon>0$ and $\alpha\in\left(
0,1\right)  ,$ there exist $r_{0},m\in%
%TCIMACRO{\U{2115} }%
%BeginExpansion
\mathbb{N}
%EndExpansion
$ satisfying%
\[
\text{ }\frac{1}{h_{r}}\sum_{k\in J_{r}}\nu_{x_{k}-x_{m}}(\varepsilon)%
 >1-\varepsilon
\]
for all $r\geq r_{0}.$\\

\textbf{Definition 3.5.} Let $I$ be an admissible ideal of $\mathbb{N}.$ Let $\left(  X,\nu,T\right)  $ be an
PNS. A sequence $x=\left(  x_{k}\right)  $ in $X$ is said to be
$I_{\theta}$-Cauchy sequence with respect to the probabilistic norm
$\nu$ if, for every $\varepsilon>0$ and
$\alpha\in\left(  0,1\right)  ,$ there exists $m\in%
%TCIMACRO{\U{2115} }%
%BeginExpansion
\mathbb{N}
%EndExpansion
$ satisfying%
\[
\left\{  r\in%
%TCIMACRO{\U{2115} }%
%BeginExpansion
\mathbb{N}
%EndExpansion
:\text{ }\frac{1}{h_{r}}\sum_{k\in J_{r}}\nu_{x_{k}-x_{m}}(\varepsilon)
 >1-\varepsilon\right\}  \in F\left(
I\right)
\]

\textbf{Definition 3.6.} Let $I$ be an admissible ideal of $\mathbb{N}.$ Let $\left(  X,\nu,T\right)  $ be an
PNS. A sequence $x=\left(  x_{k}\right)  $ in $X$ is said to
be\ $I_{\theta}^{\ast}$-Cauchy sequence with respect to the probabilistic norm $\nu$ if, there exists a set
$M=\left\{  m_{1}<m_{2}<...<m_{k}<...\right\}  \subset%
%TCIMACRO{\U{2115} }%
%BeginExpansion
\mathbb{N}
%EndExpansion
$ such that the set $M^{\imath}=\left\{  r\in%
%TCIMACRO{\U{2115} }%
%BeginExpansion
\mathbb{N}
%EndExpansion
:\text{ }m_{k}\in J_{r}\right\}  \in F\left(  I\right)  $ and the subsequence
$\left(  x_{m_{k}}\right)  $ of $x=\left(  x_{k}\right)  $ is a $\theta$-Cauchy
sequence with respect to the probabilistic norm $\nu.$\\

 The following theorem is an analogue of Theorem 3.3, so the proof omitted.\\

\textbf{Theorem 3.10. } {\it Let $I$ be an admissible ideal of $\mathbb{N}.$ Let $\left(  X,\nu,T\right)  $ be an
PNS. If a sequence $x=\left(  x_{k}\right)  $ in $X$ is $\theta$-Cauchy sequence
with respect to the probabilistic norm $\nu,$ then it is $I_{\theta}$-Cauchy sequence with respect to the
same norm.}\\

The proof of the following theorem's proof is similar to that of
Theorem 3.5.\\

\textbf{Theorem 3.11.} {\it Let $\left(  X,\nu,T\right)  $ be an
PNS. If a sequence $x=\left(  x_{k}\right)  $ in $X$ is $\theta$-Cauchy sequence
with respect to the probabilistic norm $\nu,$ then there is a subsequence of $x=\left(  x_{k}\right)  $
which is ordinary Cauchy sequence with respect to the same norm.}\\

The following theorem can be proved easily using similar techniques as in the proof of Theorem 3.6.\\

\textbf{Theorem 3.12.} {\it Let $I$ be an admissible ideal of $\mathbb{N}.$ Let $\left(  X,\nu,T\right)  $ be an PNS. If a sequence $x=\left(  x_{k}\right)  $ in $X$ is $I_{\theta}^{\ast}$-Cauchy sequence with respect to theprobabilistic norm $\nu,$ then it is $I_{\theta}$-Cauchy sequence as well.}

%\noindent {\bf Problem 1.} For further study, we suggest to investigate $\lambda$-ideal convergence for the fuzzy points. However due to the change in settings, the definitions and methods of proofs will not always be analogous to these of the present work(for example see \cite{CakalliandPratul}).

%\noindent {\bf Problem 2.} For another further study we suggest to introduce a new concept in dynamical systems using $\lambda$-ideal convergence.

%\section {conclusion}

%The present work contains not only an improvement and a generalization of the works of \cite{MursaleenMohiuddine} and \cite{Rahmat}, and 
%Section 2 of the paper \cite{MursaleenMohiuddine} as it has been presented in more general setting, i.e. in $I$-convergence which is more general than the statistical case, but also an 
%investigated some futher results in $I_{\theta}$-convergence in probabilistic normed spaces. So that one may expect it to be more useful tool in the field of metric space theory in modeling various problems occurring in many areas of science, computer science, information theory, dynamical systems, biological science, geographic information systems, population modeling, and motion planning in robotics. 
%It seems that an investigation of the present work taking "nets" instead of "sequences" could be done using the properties of "nets" instead of using the properties of "sequences". 


\begin{thebibliography}{9}    

\bibitem{AlsinaSchweizerSklar93} C. Alsina,  B. Schweizer, A. Sklar, On the definition of a probabilistic normed space. Aequationes Math.,
46, 91-98(1993)                                                                                            
\bibitem{AlsinaSchweizerSklar97} C. Alsina, B. Schweizer, A. Sklar,  Continuity properties of probabilistic norms, J. Math. Anal. Appl. 208 (1997), 446-452.
\bibitem{Buck} R.C. Buck, The measure theoretic approach to density, Amer. J. Math. 68 (1946) 560-580.
%\bibitem{CakalliandPratul} H.\c{C}akall\i, and Pratulananda Das, Fuzzy compactness via summability, Appl. Math. Lett., 22(11)(2009), 1665-1669, MR \textbf{2010k}:54006.
\bibitem{Cakalli96} H. \c{C}akalli, On statistical convergence in topological groups, Pure Appl. Math. Sci.,43(1996), 27-31.
\bibitem{Cakalli09} H. \c{C}akalli, A study on statistical convergence, Funct. Anal. Approx. Comput., 1(2)(2009), 19-24, MR2662887.
\bibitem{CakalliHazarika}  H.\c{C}akall\i, and  B. Hazarika, Ideal-quasi-Cauchy sequences, Jour. Ineq. Appl., 2012(2012) pages 11, doi:10.1186/1029-242X-2012-234
\bibitem{CasertaMaioKocinac} A. Caserta, G. Di Maio, Lj. D. R. Ko\u{c}inac, Statistical convergence in function spaces, Abstr. Appl. Anal. Vol. 2011(2011), Article ID 420419, 11 pages.
%\bibitem{CasertaKocinac} A. Caserta, Lj.D.R. Ko\u{c}inac, On statistical exhaustiveness, Appl. Math. Letters, in press.
\bibitem{ChengLinLanLiu} L. X. Cheng, G. C. Lin, Y. Y. Lan, H. Liu, Measure theory of statistical convergence, Science in China, Ser. A: Math. 51(2008), 2285-2303.
\bibitem{Choudhary09} B. Choudhary, Lacunary $I$-convergent sequences, in: Real Analysis Exchange Summer Symposium, 2009, pp. 56-57.
\bibitem{Connor} J. Connor, The statistical and strong $p$-Ces\'aro convergence of sequences, Analysis 8 (1988) 47-63.
\bibitem{Con} J.Connor, M.A. Swardson, Measures and ideals of $C^{\ast}(X) $ , Ann. N.Y. Acad.Sci.704(1993), 80-91.
\bibitem{ConstantinIstratescu} G. Constantin, I. Istratescu,  Elements of Probabilistic Analysis, Kluwer, 1989
\bibitem{Debnath} P. Debnath, Lacunary ideal convergence in intuitionistic fuzzy normed linear spaces, Comput.
Math. Appl., 63(2012),708-715.
\bibitem{Dems} K. Dems, On $I$-Cauchy sequences. Real Anal. Exchange, 30(1)(2004/2005), 123-128 
\bibitem{EsiHazarika} A. Esi, B. Hazarika, $\lambda$-ideal convergence in intuitionistic fuzzy 2-normed linear space, Jour. Intell. Fuzzy Systems, 24(4)(2013), 725-732, DOI: 10.3233/IFS-2012-0592
\bibitem{EsiHazarika00} A. Esi, B. Hazarika, Lacunary Summable Sequence Spaces of Fuzzy Numbers Defined By Ideal Convergence and an Orlicz Function, Afrika Matematika, DOI: 10.1007/s13370-012-0117-3 
\bibitem{Fast} H. Fast, Sur la convergence statistique, Colloq. Math. 2(1951) 241-244.
\bibitem{FreedmanSemberRaphael78} A. R. Freedman, J. J. Sember, M. Raphael, Some Cesaro-type summability spaces, Proc. London Math. Soc., 37(3) (1978) 508-520.
\bibitem{Fridy85} J. A. Fridy, On statistical convergence, Analysis, 5(1985) 301-313.  
\bibitem{FridyOrhan93} J. A. Fridy, C. Orhan, Lacunary statistical convergence, Pacific J. Math., 160(1)(1993), 43-51, MR 94j:40014.

\bibitem{FridyOrhanLacunary statistical summability} J. A. Fridy, C. Orhan, Lacunary statistical summability, J. Math. Anal. Appl., 173, (1993), 497-504, MR 95f :40004.
%\bibitem{Fridy93} J.A. Fridy, Statistical limit points, Proc. Amer. Math. Soc. 118 (1993) 1187-1192.
%\bibitem{Gahler65} S. Gahler, Linear 2-normietre Raume, Math. Nachr. 28(1965) 1-43.

\bibitem{GuillenLallenaSempi} B.L. Guill\'en, J.A.R. Lallena, C. Sempi, A study of boundedness in probabilistic
normed spaces. J. Math. Anal. Appl. 232 (1999), 183-196.
\bibitem{GuillenSempi} B.L. Guill\'en, C. Sempi, Probabilistic norms and convergence of random variables.
J. Math. Anal. Appl. 280 (2003), 9-16.
 %\bibitem{ConstantinIstratescu} G. Constantin, I. Istratescu,  Elements of Probabilistic Analysis, Kluwer, 1989
%\bibitem{Hazarika} B. Hazarika, Fuzzy real valued lacunary $I$-convergent sequences, Appl. Math. Letters 25 (2012) 466–470
\bibitem {Hazarika09} B. Hazarika, Lacunary $I$-convergent sequence of fuzzy real numbers, The Pacific Jour. Sci. Techno., 10(2) (2009),203-206.

\bibitem{Hazarika} B. Hazarika, Fuzzy real valued lacunary $I$-convergent sequences, Appl. Math. Letters 25 (2012) 466–470.

\bibitem{HazarikaSavas} B. Hazarika, E. Sava\c{s}, Some $I$-convergent lambda-summable difference sequence spaces of fuzzy real numbers defined by a sequence of Orlicz functions, Math. Comp. Modell.54(2011) 2986-2998.
%\bibitem{KarakusDemirciDuman} S. Karakus, K. Demirci and O. Duman, Statistical convergence on intuitionistic fuzzy normed spaces, Chaos, Solitons and Fractals, 35 (2008), 763-769.
\bibitem{HazarikaSavas123} B. Hazarika, E. Savas, Lacunary statistical convergence of double sequences and some inclusion results in $n$-normed spaces, Acta Mathematica Vietnamica, (Accepted for publications).
\bibitem {Hazarika622} B. Hazarika, Lacunary difference ideal convergent sequence spaces of fuzzy numbers, Journal of Intelligent and Fuzzy Systems,  DOI: 10.3233/IFS-2012-0622.
\bibitem {Hazarika12} B. Hazarika, On generalized difference ideal convergence in random 2-normed spaces, Filomat, 26(6) (2012), 1265-1274.
\bibitem{Hazarika769} B. Hazarika, On $\sigma$-uniform density and ideal convergent sequences of fuzzy real numbers, Journal of Intelligent and Fuzzy Systems, doi 10.3233/IFS-130769.
\bibitem{HazarikaKimarGuillen} B. Hazarika, V. Kumar, B. L. Guill\'en, Generalized ideal convergence in intuitionistic fuzzy normed linear spaces, Filomat (Accepted for publications).
\bibitem{HazarikaOn ideal convergence in topological groups} B. Hazarika, On ideal convergence in topological groups, Scientia Magna,7(4)(2011), 80-86.
\bibitem{Karakus} S. Karakus, Statistical convergence  on probabilistic normed spaces. Math. Comm. 12
(2007), 11-23.
\bibitem{KlementMesiarPap}  E. P. Klement, R. Mesiar, E. Pap,, Triangular Norms, Kluwer, Dordrecht, 2000
\bibitem{KostyrkoSalatWilczynski} P. Kostyrko, T. \u{S}al\'at, and W. Wilczy\'nski, $I$-convergence, Real Anal. Exchange \textbf{26}, 2, (2000-2001),669-686, MR \textbf{2002e}:54002.
 %\bibitem{ AlsinaSchweizerSklar} C. Alsina,  B. Schweizer, A. Sklar, On the definition of a probabilistic normed space. Aequationes Math.,46, 91–98(1993)

\bibitem{KostyrkoMacajSalatSleziak} P. Kostyrko, M. Macaj, T. \u{S}alat, M. Sleziak, $I$-convergence and Extremal $I$-limit 
        points, Math. Slovaca 2005; 55; 443-64.
 
 
  
%\bibitem{Sherstnev} A. N. \u{S}herstnev, Random normed spaces: Problems of completeness. Kazan Gos. Univ. Ucen. Zap., 122, 3–20 (1962)       
        %\bibitem{DasKostyrkoWilczynskiMalik} P. Das, P. Kostyrko, W. Wilczy\'nski, P. Malik, $I$ and $I^{\ast}$-convergence of double sequences, Math.Slovaca, 58 (5) (2008), 605-620.
      
        
%\bibitem{KumarKumar} V. Kumar, K. Kumar, On the ideal convergence of sequences in intuitionistic fuzzy normed spaces, Selcuk J. Math. 2009; 10(2); 27-41.
\bibitem{KumarKumar08} K. Kumar, V. Kumar, On the $I$ and $I^{\ast}$-Cauchy sequences in probabilistic normed spaces. Mathematical
Sciences, 2(1), 47-58 (2008)
%\bibitem{Kumar} V. Kumar, On $I$ and $I^{\ast}$-convergence of double sequences, Math. Commun. 12 (2007), 171.181.
%\bibitem{KumarKumar08} K. Kumar, V. Kumar, On the $I$ and $I^{\ast}$-convergence of sequences in fuzzy normed Spaces, Advances in Fuzzy Sets and Systems, 2008, 3(3), 341-365.

%\bibitem{KumarMursaleen} V. Kumar, M. Mursaleen, On $(\lambda, \mu)$-statistical convergence of double sequences on intuitionistic fuzzy normed spaces, Filomat 25:2 (2011), 109-120. DOI: 10.2298/FIL1102109K.
\bibitem{KumarGuillen} V. Kumar, B. L. Guill\'en, On Ideal Convergence of
Double Sequences in Probabilistic Normed Spaces, Acta Math. Sinica, English
Series, Published online: February 21, 2012, DOI: 10.1007/s10114-012-9321-1

\bibitem{LahiriDas} B. K. Lahiri, P. Das, $I$ and $I^{\ast}$-convergence in topological spaces, Math. Bohemica, 130 (2005), 153-160.
\bibitem{LafuerzaLallenaSempi} B. Lafuerza-Guill\'en,  J. A.Rodr\'iguez-Lallena, C. Sempi,  Some classes of probabilistic normed spaces.
Rend. Mat., 17(1997),237-252 

%\bibitem {Leindler} L.Leindler, \"{U}ber die de la Vallee-Pousinsche Summeierbarkeit allgemeiner Orthogonalreihen, Acta Math.Acad.Sci.Hungar., 16(1965), 375-387.
\bibitem{Li00} J. Li, Lacunary statistical convergence and inclusion properties between lacunary methods, Internat. J. Math.  Math. Sci. 23(3) (2000), 175-180, S0161171200001964.
\bibitem{Maddox88} I. J. Maddox, Statistical convergence in a locally convex spaces, Math. Proc. Cambridge Philos.
Soc., 104(1)(1988), 141-145.
\bibitem{MaioKocinac08} G. Di. Maio,  Lj.D.R. Ko\u{c}inac, Statistical convergence in topology, Topology Appl. 156,
(2008), 28-45.

\bibitem{Menger} K. Menger, Statistical metrics. Proc. Nat. Acad. Sci. USA, 28(1942) 535-537 
\bibitem{Miller95} H. I. Miller, A measure theoretical subsequence characterization of statistical convergence,
Trans. Amer. Math. Soc., 347(5)(1995), 1811-1819.
\bibitem{MitrinovicSandorCrstici} D.S. Mitrinovi\'c, J. Sandor, B. Crstici, Handbook of Number Theory, Kluwer Acad. Publ., Dordrecht, Boston, London, 1996.
%\bibitem{MursaleenAlotaibiStatisticalsummabilityandapproximationbydelaValleePousinmean} M. Mursaleen and A. Alotaibi, Statistical summability and approximation by de la Vallee-Pousin mean, Applied Math. Letters, 24 (2011) 320-324.

%\bibitem{MursaleenEdelyGeneralizedstatisticalconvergence} M. Mursaleen and Osama H.H.Edely, Generalized statistical convergence, Information Sciences , 162(2004), 287-294.

%\bibitem{MursaleenEdelyOntheinvariantmeanandstatisticalconvergence} M. Mursaleen and Osama H.H. Edely, On the invariant mean and statistical convergence, Appl. Math. Letters,22(2009) 1700-1704.

%\bibitem{MohiuddineLohani} S.A. Mohiuddine and Q.M. Danish Lohani, On generalized statistical convergence in intuitionistic fuzzy normed space, Chaos, Solitons and Fractals, 42 (2009) 1731-1737.

%\bibitem{MursaleenMohiuddineOnlacunarystatisticalconvergencewithrespecttotheintuitionisticfuzzynormedspace} M. Mursaleen and S.A. Mohiuddine, On lacunary statistical convergence with respect to the intuitionistic fuzzy normed space, Jour. Comput. Appl. Math., 233(2)(2009) 142-149.

%\bibitem{MursaleenMohiuddineEdely} M. Mursaleen, S. A. Mohiuddine, Osama H.H.Edely, On ideal convergence of double sequences in intuitionistic fuzzy normed spaces, Comput. Math. Appl. 59 (2010) 603-611.
%\bibitem{Mursaleen} M. Mursaleen, $\lambda$-statistical convergence, Mathematica Slovaca, 50(1)(2000),111-115.
\bibitem{MursaleenMohiuddine} M. Mursaleen, S. A. Mohiuddine, On ideal convergence in probabilistic normed spaces, Math. Slovaca, 62(1)(2012), 49-62. 
%\bibitem{MursaleenLohani08} M.Mursaleen and Q.M.Danish Lohani, Intuitionistic fuzzy 2-normed space and some related concepts, Chaos, Solitons and Fractals, (2008), doi:10.1016/j.chaos.2008.11.006.

%\bibitem{SahinerGurdalSaltanGunawan} A. Sahiner, M.G\"urdal, S. Saltan and H. Gunawan, Ideal Convergence in 2-normed spaces,Taiwanese J. Math., 11(5) 2007; 1477 ¡ 1484:
\bibitem{Rahmat} M. R. S. Rahmat, Ideal Convergence on Probabilistic Normed Spaces, Inter. Jour. Stat. Econ.,3(9)(2009), 67-75
\bibitem{Salat80} T. \u{S}al\'at, On statistical convergence of real numbers, Math. Slovaca, 30(1980), 139-150.

%\bibitem{LafuerzaLallena Sempi} B. Lafuerza-Guill\'en,  J. A.Rodr\'iguez-Lallena, C. Sempi,  Some classes of probabilistic normed spaces.Rend. Mat., 17, 237–252 (1997)
\bibitem{SalatTripathyZiman04} T. \u{S}al\'at, B. C. Tripathy, M. Ziman, On some properties of $I$-convergence, Tatra 
        Mt. Math. Publ. 2004; 28; 279-86.

%\bibitem{SalatTripathyZiman05} T. \u{S}al\'at, B. C. Tripathy, M. Ziman, On $I$-convergence field, Italian J. Pure and        Appl. Math. 2005; 17; 45-54. 
       %\bibitem{SalatTijdeman} T. \u{S}al\'at, R. Tijdeman, On statistically convergent sequences of real numbers, Math. Slovaca 30 (1980) 139–150. 
       %\bibitem{Saminger Sempi08} S. Saminger-Platz, C. Sempi,  A primer on triangle functions I. Aequationes Math., 76, 201-240 (2008)
%\bibitem{Saminger Sempi10}S. Saminger-Platz, C.Sempi, : A primer on triangle functions II. Aequationes Math., 80(3), 239-268 (2010)
%\bibitem {SavasDas} E. Savas, P. Das, A generalized statistical convergence via ideals, Applied Math. Letters, 24(2011), 826-830.

%\bibitem{Schoenberg59} I. J. Schoenberg, The integrability of certain functions and related summability methods, Amer. Math. Monthly 66(1959) 361-375.

\bibitem{SchweizerSklar60} B. Schweizer, A. Sklar, Statistical metric spaces, Pacific J. Math. 10 (1960),
313-334.

\bibitem{SchweizerSklar83} B. Schweizer, A. Sklar, Probabilistic Metric Spaces, North Holland, New York- Amsterdam-Oxford, 1983.
\bibitem{Sherstnev} A. N. \u{S}herstnev, Random normed spaces: Problems of completeness. Kazan Gos. Univ. Ucen. Zap., 122,
3-20 (1962)
\bibitem{Steinhaus} H. Steinhaus, Sur la convergence ordinaire et la convergence asymptotique, Colloq. Math. 2 (1951) 73-74.  
        
%\bibitem {1} P. Kostyrko, T.Salat and W. Wilczynski, I-convergence, Real Analysis Exchange, 26(2)(2000/2001), 669-686.

%\bibitem {2} T. Salat, B. C. Tripathy and M. Ziman, On I-convergence field, Italian J. Pure and Appl. Math., 17(2005), 45-54.


%\bibitem{Saminger Sempi08} S. Saminger-Platz, C. Sempi,  A primer on triangle functions I. Aequationes Math., 76, 201–240 (2008)
%\bibitem{Saminger Sempi10}S. Saminger-Platz, C.Sempi, : A primer on triangle functions II. Aequationes Math., 80(3), 239–268 (2010)





%\bibitem {4} H.Fast, Sur la convergence statistique, Colloq.Math., 2(1951), 241-244.


%\bibitem{MursaleenMohiuddineEdely} M. Mursaleen, S. A. Mohiuddine, On ideal convergence in probabilistic normed spaces, Math. Slovaca, 62(1)(2012), 49-62. 
 

%\bibitem{Rahmat} M. R. S. Rahmat, Ideal Convergence on Probabilistic Normed Spaces, Inter. Jour. Stat. Econ.,3(9)(2009), 67-75

%\bibitem {6 }B. Schweizer and Sklar, A, Statistical metric spaces, Pacific J.Math., 10(1960), 313-334.


%\bibitem{Gahler63} S. G\"{a}hler, 2-metrische Raume and ihre topologische Struktur, Math. Nachr. 26(1963) 115-148.



%\bibitem {7}S. Gahler, Lineare 2-normietre Raume, Math.Nachr, 28(1965), 1-43.


%\bibitem{TripathyTripathy} B. K. Tripathy, B.C. Tripathy, On $I$-convergent double sequences, Soochow J. Math., 31(4) (2005) 549-560. 
       
\bibitem{TripathyHazarika} B. C. Tripathy, B. Hazarika,  $I$-monotonic and $I$-convergent sequences, 
        Kyungpook Math. J. 51(2011), 233-239, DOI 10.5666/KMJ.2011.51.2.233.
        
        %\bibitem{Choudhary09} B. Choudhary, Lacunary $I$-convergent sequences, in: Real Analysis Exchange Summer Symposium, 2009, pp. 56-57.

\bibitem{TripathyHazarikaChoudhary} B. C. Tripathy, B. Hazarika, B. Choudhary, Lacunary $I$-convergent sequences, Kyungpook Math. J. 52(4)(2012) 473-482.
\bibitem{TripathyHazarika09} B. C. Tripathy and B. Hazarika, Paranorm $I$-convergent sequence
spaces, \textit{Math. Slovaca,} \textbf{59(4)} (2009) 485-494.
\bibitem{TripathyHazarika11} B. C. Tripathy, B. Hazarika, Some $I$-convergent sequence spaces
defined by Orlicz functions, \textit{Acta Math. Appl. Sinica,} \textbf{27(1)}
(2011) 149-154.
\bibitem{YamanciGurdal} U. Yamanc\i, M. G\"{u}rdal, On lacunary ideal convergence in random $n$-normed space, Journal of Mathematics, Vol. 2013(2013), Article ID 868457, 8 pages .
  
\end{thebibliography}
\end{document}